\documentclass[11pt]{article}
\usepackage[matrix,arrow,curve,cmtip]{xy}
\usepackage{amssymb}
\usepackage{latexsym}
\usepackage{theorem}
\usepackage[a4paper,width=160mm,height=230mm]{geometry}
\usepackage{titlesec}
\usepackage{graphicx}

\titleformat{\section}[hang]%
{\bf\Large}{\thesection.}{1ex}{}%

\titleformat{\subsection}[hang]%
{\bfseries\normalsize}{\thesubsection.}{1ex}{}

\def\to{\mbox{$\xymatrix@1@C=5mm{\ar@{->}[r]&}$}}
\def\halfcirc{\begin{picture}(0,0)\put(0,3){\oval(4,4)[l]}\end{picture}}
\def\incl{\mbox{$\xymatrix@1@C=5mm{\ar@{->}[r]|<{\halfcirc}&}$}}
\def\tto{\mbox{$\xymatrix@1@C=5mm{\ar@{=>}[r]&}$}}
\def\distsign{\begin{picture}(0,0)\put(0,0){\circle{4}}\end{picture}}
\def\dist{\mbox{$\xymatrix@1@C=5mm{\ar@{->}[r]|{\distsign}&}$}}
\def\spansign{\begin{picture}(0,0)\put(0,-3){\line(0,1){6}}\end{picture}}
\def\span{\mbox{$\xymatrix@1@C=5mm{\ar@{->}[r]|{\spansign}&}$}}
\def\criblesign{\begin{picture}(0,0)\put(-1,-3){\line(0,1){6}}\put(1,-3){\line(0,1){6}}\end{picture}}
\def\crible{\mbox{$\xymatrix@1@C=5mm{\ar@{->}[r]|{\criblesign}&}$}}

\def\inlineadj#1#2{\mbox{$\xymatrix@C=15mm{\ar@{}[r]|{\bot}\ar@<1mm>@/^2mm/[r]^{{#1}} &\ar@<1mm>@/^2mm/[l]^{{#2}}}$}}

\newtheorem{theorem}{Theorem}[section]
\newtheorem{lemma}[theorem]{Lemma}
\newtheorem{definition}[theorem]{Definition} 
\newtheorem{proposition}[theorem]{Proposition}

{\theorembodyfont{\upshape}\newtheorem{example}[theorem]{Example}}
\newcommand{\proof}{\noindent {\em Proof\ }: }
\def\endofproof{$\mbox{ }\hfill\Box$\par\vspace{1.8mm}\noindent}

\def\Rel{{\sf Rel}}
\def\o{^{\sf o}}
\def\bigmid{~\Big|~}

\def\si{_{\sf si}}

\def\Matr{{\sf Matr}}

\def\:{\colon}

\def\2{{\bf 2}}
\def\Set{{\sf Set}}
\def\op{^{\sf op}}
\def\dom{{\sf dom}}
\def\cod{{\sf cod}}
\def\Sup{{\sf Sup}}
\def\Dist{{\sf Dist}}
\def\Map{{\sf Map}}

\def\id{{\sf id}}

\def\C{{\cal C}}

\def\G{{\cal G}}

\def\Q{{\cal Q}}
\def\R{{\cal R}}

\def\bbA{\mathbb{A}}
\def\bbB{\mathbb{B}}
\def\bbC{\mathbb{C}}

\def\bbR{\mathbb{R}}

\def\tensor{\otimes}
\def\<{\langle}
\def\>{\rangle}

\def\eqref#1{(\ref{#1})}
\def\inv{^{-1}}

\def\obis{^{\sf x}}

\def\down{\downarrow\!}
\def\otspam{\begin{picture}(16,10)\put(2,5){\rotatebox{180}{$\mapsto$}}\end{picture}}

\title{Elementary characterisation of small quantaloids of closed cribles}
\author{Hans Heymans\footnote{Department of Mathematics and Computer Science, University of Antwerp, Middelheimlaan 1, 2020 Antwerpen, Belgium, {\tt hans.heymans@ua.ac.be}}\ \ and Isar Stubbe\footnote{Laboratoire de Math\'ematiques Pures et Appliqu\'ees, Universit\'e du Littoral-C\^ote d'Opale, 50 rue F. Buisson, 62228 Calais, France, {\tt isar.stubbe@lmpa.univ-littoral.fr}}}
\date{March 17, 2011}

\begin{document}

\maketitle

\begin{abstract}
Each small site $(\C,J)$ determines a small quantaloid of closed cribles $\R(\C,J)$. We prove that a small quantaloid $\Q$ is equivalent to $\R(\C,J)$ for some small site $(\C,J)$ if and only if there exists a (necessarily subcanonical) Grothendieck topology $J$ on the category $\Map(\Q)$ of left adjoints in $\Q$ such that $\Q\cong\R(\Map(\Q),J)$, if and only if $\Q$ is locally localic, map-discrete, weakly tabular and weakly modular. If moreover coreflexives split in $\Q$, then the topology $J$ on $\Map(\Q)$ is the canonical topology.
\end{abstract}

\section{Introduction}\label{A}

A quantaloid $\Q$ is, by definition, a category enriched in the symmetric monoidal closed category $\Sup$ of complete lattices and supremum-preserving functions [Rosenthal, 1996]. Viewing $\Q$ as a bicategory, it is natural to study categories, functors and distributors enriched in $\Q$ [B\'{e}nabou, 1967; Street, 1983; Stubbe, 2005a]. A major application of quantaloid-enriched category theory was discovered by B. Walters, and published in this journal: in [1982], he proved that the topos of sheaves on a site $(\C,J)$ is equivalent to the category of symmetric and Cauchy complete categories enriched in the {\em small quantaloid of closed cribles} $\R(\C,J)$ constructed from the given site.

Given the importance of the {\em construction} of the quantaloid of closed cribles $\R(\C,J)$ from a small site $(\C,J)$, we provide in this paper an elementary {\em axiomatisation} of this notion. Precisely, we prove that a small quantaloid $\Q$ is equivalent to $\R(\C,J)$ for some small site $(\C,J)$ if and only if there exists a Grothendieck topology $J$ on the category $\Map(\Q)$ of left adjoints in $\Q$ such that $\Q\cong\R(\Map(\Q),J)$, if and only if $\Q$ is locally localic, map-discrete, weakly tabular and weakly modular. (The latter two notions seem to be new, and inherited their name from the stronger notions of tabularity and modularity introduced in [Freyd and Scedrov, 1990].) The Grothendieck topology $J$ on $\Map(\Q)$ is always subcanonical, and if coreflexives split in $\Q$, then $J$ is the canonical topology.

This result thus spells out how two, at first sight quite different, generalisations of locales, namely Grothendieck topologies on the one hand, and quantaloids on the other, relate: the former can be understood to form an axiomatically described subclass of the latter. It is hoped that this axiomatisation helps to clarify the role that quantaloids may play in the search for a good notion of ``non-commutative topology'', to be used ultimately in suitable generalisations of sheaf theory (see e.g.\ [Borceux and Van den Bossche, 1986; Mulvey and Nawaz, 1995; Ambler and Verity, 1996; H\"{o}hle, 1998; Gylys, 2001; Garraway, 2005; Stubbe, 2005b]).

\section{Small quantaloids of closed cribles}\label{B}

To begin, we recall a construction due to [Walters, 1982]. If $\C$ is a small category, then the  quantaloid $\R(\C)$ of {\em cribles} in $\C$ is the full sub-quantaloid of $\Rel(\Set^{\C\op})$ whose objects are the representable presheaves. It is useful to have an explicit description. We write a {\em span} in $\C$ as $(f,g)\:X\span Y$, and intend it to be a pair of arrows with $\dom(f)=\dom(g)$, $\cod(f)=Y$ and $\cod(g)=X$:
$$\xymatrix@C=3ex{
 & \cdot\ar[dl]_g\ar[dr]^f & \\
X\ar[rr]|{\spansign}_{(f,g)} & & Y}$$
(Many would consider such a pair to be a span in the opposite direction, but the reason for our notational convention for domain and codomain will become clear when we compose cribles.) A crible $R\:X\crible Y$ is then a set of spans $X\span Y$ such that for any $(f,g)\in R$ and any $h\in\C$ with $\cod(h)=\dom(f)$, also $(f\circ h,g\circ h)\in R$. Composition in $\R(\C)$ goes as follows: for $R\:X\crible Y$ and $S\:Y\crible Z$ the elements of $S\circ R\:X\crible Z$ are those $(f,g)$ for which there exists a morphism $t\in\C$ such that $(f,t)\in S$ and $(t,g)\in R$. The identity crible $\id_X\:X\crible X$ is the set $\{(f,f)\mid\cod(f)=X\}$. (Here we need $\C$ to be small: otherwise $\id_X$ isn't necessarily a {\em set}.) The supremum of a set of cribles from $X$ to $Y$ is simply their set-theoretic union. 

Next recall from [Mac Lane and Moerdijk, 1992] that a {\em Grothendieck topology $J$} on a small category $\C$ is a function, assigning to every object $C$ a set $J(C)$ of sieves on $C$, that satisfies three conditions:
\begin{itemize}
\setlength{\itemsep}{0ex}
\item $\top_C:=\{f\in\C\mid\cod(f)=C\}\in J(C)$,
\item if $S\in J(C)$ then $f^*(S):=\{g\in\C\mid f\circ g\in S\}\in J(D)$ for any $f\:D\to C$ in $\C$,
\item if $S\in J(C)$ and $T$ is a sieve on $C$ such that $s^*(T)\in J(\dom(s))$ for all $s\in S$, then $T\in J(C)$ too.
\end{itemize}
An element of $J(C)$ is a {\em covering sieve} on $C$; the couple $(\C,J)$ is a {\em small site}.

A {\em nucleus} $j$ on a quantaloid $\Q$ is a lax functor $j\:\Q\to\Q$ which is the identity on objects and such that each $j\:\Q(X,Y)\to\Q(X,Y)$ is a closure operator; it is {\em locally left exact} if it preserves finite infima of arrows. Grothendieck topologies $J$ on $\C$ are in bijective correspondence with locally left exact nuclei $j$ on $\R(\C)$ [Betti and Carboni, 1983; Rosenthal, 1996], as follows: For a Grothendieck topology $J$ on $\C$,
let $j\:\R(\C)\to\R(C)$ send a crible $R\:C\crible D$ to
$$j(R):=\left\{(f,g)\:C\span D\bigmid\exists S\in J(\dom(f)):\forall s\in S,\ (g\circ s,f\circ s)\in R\right\}.$$
Conversely, if $j\:\R(\C)\to\R(\C)$ is a locally left exact nucleus, then put
$$J(C):=\left\{S\mbox{ is a sieve on }C\bigmid\id_C\leq j(\{(s,s)\mid s\in S\})\right\}.$$
If $j\:\Q\to\Q$ is a nucleus on a quantaloid, then there is a (``quotient'') quantaloid $\Q_j$ of {\em $j$-closed morphisms}, i.e.\ those $f\in\Q$ for which $j(f)=f$: the composition is $j(g\circ f)$, the identity on an object $X$ is $j(1_X)$, and the supremum of a family $(f_i)_{i\in I}$ is $j(\bigvee_if_i)$. For a small site $(\C,J)$ we write $\R(\C,J)$ for the quantaloid $\R(\C)_j$ with $j$ the nucleus determined by the topology $J$. 
\begin{definition}[Walters, 1982]\label{1}
A {\em small quantaloid of closed cribles} is a small quantaloid which is equivalent to $\R(\C,J)$ for some small site $(\C,J)$. 
\end{definition}
To be precise, Walters [1982] called this a `bicategory of relations', wrote it as $\Rel(\C,J)$, and called its arrows `relations'. However, to avoid confusion with the `bicategories of relations' that [Carboni and Walters, 1987] and others have since then worked on, we prefer to stick closer to the actual construction and speak of a `quantaloid of closed cribles'.

In the remainder of this paper we develop an axiomatic description of the class of small quantaloids of closed cribles, {\em purely in terms of composition and local suprema/infima}. We start by preparing the ground in the next section, in which we indicate several key properties of such quantaloids.

\section{Weak tabularity, weak modularity}\label{C}

Recall that an {\em involution} on a quantaloid $\Q$ is a $\Sup$-functor $(-)\o\:\Q\op\to\Q$ which is the identity on objects and satisfies $f^{\sf oo}=f$ for any morphism $f$ in $\Q$; the pair $(\Q,(-)\o)$ is then said to form an {\em involutive quantaloid}, but most of the time we do not explicitly mention the functor $(-)\o$ and simply speak of `an involutive quantaloid $\Q$'. If a morphism $f\:X\to Y$ in a quantaloid $\Q$ is a left adjoint, then we write its right adjoint as $f^*\:Y\to X$. To avoid overly bracketed expressions we often write $gf$ instead of $g\circ f$ for the composition of two morphisms $f\:X\to Y$ and $g\:Y\to Z$ in $\Q$. The next definition gathers some technical conditions which will show up in the main theorem in the next section.
\begin{definition}\label{2}
A quantaloid $\Q$ is:
\begin{enumerate}
\item\label{2.1} {\em locally localic} if, for all objects $X$ and $Y$, $\Q(X,Y)$ is a locale,
\item\label{2.2} {\em map-discrete} if, for any left adjoints $f\:X\to Y$ and $g\:X\to Y$ in $\Q$, $f\leq g$ implies $f=g$,
\item\label{2.3} {\em weakly tabular} if, for every $q\:X\to Y$ in $\Q$, 
$$q=\bigvee\left\{fg^*\bigmid (f,g)\:X\span Y\mbox{ is a span of left adjoints such that }fg^*\leq q\right\},$$
\item\label{2.4} {\em map-tabular} if for every $q\:X\to Y$ in $\Q$ there is a span $(f,g)\:X\span Y$ of left adjoints in $\Q$ such that $fg^*=q$ and $f^*f\wedge g^*g=1_{\dom(f)}$,
\item\label{2.5} {\em weakly modular} if, for every pair of spans of left adjoints in $\Q$, say $(f,g)\:X\span Y$ and $(m,n)\:X\span Y$, we have $fg^*\wedge mn^*\leq f(g^*n\wedge f^*m)n^*$,
\item\label{2.6} {\em tabular} if it is involutive and if for every $q\:X\to Y$ in $\Q$ there exists a span $(f,g)\:X\span Y$  of left adjoints in $\Q$ such that $fg\o=q$ and $f\o f\wedge g\o g=1_{\dom(f)}$,
\item\label{2.7} {\em modular} if it is involutive and if for any $f\:X\to Y$, $g\:Y\to Z$ and $h\:X\to Z$ in $\Q$ we have $gf\wedge h\leq g(f\wedge g\o h)$ (or equivalently, $gf\wedge h\leq (g\wedge hf\o)f$).
\end{enumerate} 
\end{definition}
The notions of modularity and tabularity are due to [Freyd and Scedrov, 1990]; we believe that weak modularity, weak tabularity and map-tabularity are new notions. Note that conditions \ref{2.1}--\ref{2.5} in the definition above make sense in {\em any} quantaloid, whereas conditions \ref{2.6}--\ref{2.7} are only defined for an {\em involutive} quantaloid. In the next two lemmas we record some straightforward implications.
\begin{lemma}\label{3.01}
If a quantaloid $\Q$ is map-tabular then it is also weakly tabular.
\end{lemma}
\begin{lemma}\label{3}
If $\Q$ is a modular quantaloid then:
\begin{enumerate}
\item if $f\:A\to B$ is a left adjoint morphism in $\Q$ then $f^*=f\o$,
\item $\Q$ is map-discrete,
\item $\Q$ is weakly modular,
\item $\Q$ is tabular if and only if it is map-tabular.
\end{enumerate}
\end{lemma}
Now we can easily point out our main example: 
\begin{example}\label{3.1}
For any small category $\C$, the quantaloid $\R(\C)$ of cribles in $\C$ is an involutive quantaloid: the involute $R\o\:D\crible C$ of a crible $R\:C\crible D$ is obtained by reversing the spans in $R$. It is easy to see that $\R(\C)$ is locally localic and modular, and by Lemma \ref{3} it is thus also map-discrete and weakly modular. Furthermore, it is weakly tabular: if we write $\<f,g\>\:C\crible D$ for the crible generated by a span $(f,g)\:C\span D$ (in the obvious way), then it is straightfoward to check that, given a crible $R\:C\crible D$, we may write
$$R=\bigcup\Big\{\<f,1\>\circ\<1,g\>\bigmid (f,g)\in R\Big\},$$
where $\<f,1\>$ is a left adjoint, and $\<1,g\>$ a right adjoint, in $\R(\C)$.

If $J$ is a Grothendieck topology on $\C$, then $\R(\C,J)$ too is involutive, because the corresponding locally left exact nucleus $j\:\R(\C)\to\R(\C)$ preserves the involution. Moreover, the involutive quantaloid $\R(\C,J)$ is locally localic and modular, because $\R(\C)$ is so and $j$ preserves these properties; thus, $\R(\C,J)$ is also map-discrete and weakly modular. Moreover, $\R(\C,J)$ is weakly tabular, again because $\R(\C)$ is so and $j$ preserves this property.
\end{example}

In the rest of this section, we relate the notions summed up in Definition \ref{2}; strictly speaking, none of these results are needed for the proof of our main theorem in the next section, but they are interesting in their own right. We start with a less straightforward relation between map-tabularity and weak tabularity in the next proposition, making use of the quantaloid $\Dist(\Q)$ of $\Q$-enriched categories and distributors (= modules = profunctors) between them. We typically write $\Phi\:\bbA\dist\bbB$ for an arrow in $\Dist(\Q)$ (whose elements are $\Q$-arrows $\Phi(b,a)\:ta\to tb$), whereas the composition with another $\Psi\:\bbB\dist\bbC$ is written as $\Psi\tensor\Phi\:\bbA\dist\bbC$ (and has elements $(\Psi\tensor\Phi)(c,b)=\bigvee_{b\in\bbB}\Psi(c,b)\circ\Phi(b,a)$). We refer to [Stubbe, 2005a] for more details and for historically relevant references.
\begin{proposition}\label{7}
A small quantaloid $\Q$ is weakly tabular if and only if $\Dist(\Q)$ is map-tabular.
\end{proposition} 
\proof
Fist suppose that $\Dist(\Q)$ is map-tabular. A $\Q$-arrow $q\:X\to Y$ may be viewed as a distributor between one-object $\Q$-categories with identity homs: $(q)\:*_X\dist *_Y$. Thus, there exist left adjoint distributors $\alpha\:\bbA\dist*_Y$ and $\beta\:\bbA\dist*_X$ satisfying in particular $\alpha\tensor\beta^*=(q)$. Spelled out, this means that 
$$q=\bigvee\Big\{\alpha(x)\circ\beta^*(x)\bigmid x\in\bbA_0\Big\}.$$ 
But it is easily seen that $\alpha(x)\dashv\alpha^*(x)$ and $\beta(x)\dashv\beta^*(x)$ for all $x\in\bbA_0$. Thus each pair $(\alpha(x),\beta(x))$ is a span of left adjoints in $\Q$, satisfying $\alpha(x)\circ\beta^*(x)\leq q$, and the above equation implies that $\Q$ is weakly tabular.

Conversely, supposing that $\Q$ is weakly tabular, we seek, for any given distributor $\Phi\:\bbB\dist\bbA$ between $\Q$-categories, left adjoint distributors $\Sigma\:\bbR\dist\bbA$ and $\Theta\:\bbR\dist\bbB$ such that $\Sigma\tensor\Theta^*=\Phi$ and $\Sigma^*\tensor\Sigma\ \wedge\ \Theta^*\tensor\Theta=\bbR$. We may suppose for convenience that $\bbA$ and $\bbB$ are Cauchy complete (because in $\Dist(\Q)$ every $\Q$-category is isomorphic to its Cauchy completion), so that any left adjoint distributor into $\bbA$ or $\bbB$ is necessarily representable. Thus our problem becomes: to find functors $S\:\bbR\to\bbB$ and $T\:\bbR\to\bbA$ such that $\bbB(-,S-)\tensor\bbA(T-,-)=\Phi$ and $\bbB(S-,S-)\wedge\bbA(T-,T-)=\bbR$. Thereto, we define the $\Q$-category $\bbR$ to be the full subcategory of $\bbA\times\bbB$ whose objects are those $(a,b)\in\bbA\times\bbB$ for which $1_{ta}\leq\Phi(a,b)$. Explicitly, $\bbR$ is given by: 
\begin{itemize}
\item objects: $\bbR_0=\{(a,b)\in\bbA_0\times\bbB_0\mid ta=tb\mbox{ and }1_{ta}\leq\Phi(a,b)\}$ with types $t(a,b)=ta=tb$,
\item hom-arrows: $\bbR((a',b'),(a,b))=\bbA(a',a)\wedge\bbB(b',b)$.
\end{itemize}
Naturally, we let $T$ (resp.\ $S$) be the composition of the inclusion $\bbR\hookrightarrow\bbA\times\bbB$ with the projection of $\bbA\times\bbB$ onto $\bbA$ (resp.\ onto $\bbB$). By construction we then have $\bbB(S-,S-)\wedge\bbA(T-,T-)=\bbR$; and a computation shows furthermore that, for any $a\in\bbA_0$ and $b\in\bbB_0$,
\begin{eqnarray*}
\bbA(a,T-)\circ\bbB(S-,b)
 & = & \bigvee_{(x,y)\in\bbR_0}\bbA(a,T(x,y))\circ\bbB(S(x,y),b) \\
 & = & \bigvee_{(x,y)\in\bbR_0}\bbA(a,x)\circ\bbB(y,b) \\
 & \leq & \bigvee_{(x,y)\in\bbR_0}\bbA(a,x)\circ\Phi(x,y)\circ\bbB(y,b) \\
 & \leq & \bigvee_{x\in\bbA}\bigvee_{y\in\bbB}\bbA(a,x)\circ\Phi(x,y)\circ\bbB(y,b) \\
 & \leq & \Phi(a,b)
\end{eqnarray*}
(using that $\Phi(x,y)\geq 1_{tx}$ to pass from the second line to the third). It remains to prove that $\Phi(a,b)\leq\bbA(a,T-)\tensor\bbB(S-,b)$ holds too. By weak tabularity of $\Q$, it suffices to show that, for any span $(f,g)\:tb\span ta$ of left adjoints in $\Q$, 
$$f\circ g^*\leq\Phi(a,b)\ \Longrightarrow \ f\circ g^*\leq\bbA(a,T-)\circ\bbB(S-,b).$$ 
Because we assumed that $\bbA$ is Cauchy complete, we can consider the tensor\footnote{By definition, the tensor $a\tensor f$ of an object $a\in\bbA$ and a morphism $f\:X\to ta$ in $\Q$, is the colimit of the functor $*_{ta}\to\bbA\:*\mapsto a$ weighted by the distributor $(f)\:*_X\dist*_{ta}$. The dual notion is cotensor; we write $\<g,a\>$ for the cotensor of an object $a\in\bbA$ with a morphism $g\:ta\to Y$ in $\Q$. If $f\dashv f^*$ in $\Q$, then $(f)$ is then a left adjoint in $\Dist(\Q)$, so a Cauchy complete $\Q$-category $\bbA$ necessarily has all tensors $a\tensor f$ with left adjoint $f$. Moreover, in this case, $a\tensor f=\<f^*,a\>$.}
$a\tensor f\in\bbA$ of the object $a\in\bbA_0$ with the left adjoint morphism $f$ in $\Q$; reckoning that $f\dashv f^*$ in $\Q$ we have moreover that the tensor $a\tensor f$ equals the cotensor $\<f^*,a\>$ of $a$ with $f^*$. A straightforward computation with the universal property of (co)tensors shows that
$$\bbA(x,\<f^*,a\>)=\bbA(x,a)\circ f\mbox{ \ and \ }\bbA(a\tensor f,y)=f^*\circ\bbA(a,y)$$
for all $x,y\in\bbA$. Similar calculations can be made for the (co)tensor $b\tensor g=\<g^*,b\>$ in $\bbB$. From this it follows easily that 
$$\Phi(a\tensor f,b\tensor g)=\bbA(a\tensor f,-)\tensor\Phi\tensor\bbB(-,\<g^*,b\>)=f^*\circ\Phi(a,b)\circ g$$
from which we can deduce that 
$$f\circ g^*\leq\Phi(a,b)\iff 1_X\leq\Phi(a\tensor f,b\tensor g)\iff(a\tensor f,b\tensor g)\in\bbR.$$ 
But this in turn implies that
\begin{eqnarray*}
\bbA(a,T-)\tensor\bbB(S-,b)
 & = & \bigvee_{(x,y)\in\bbR}\bbA(a,T(x,y))\circ\bbB(S(x,y),b) \\
 & \geq & \bbA(a,a\tensor f)\circ\bbB(b\tensor g,b) \\
 & = & \bbA(a,a)\circ f\circ g^*\circ\bbB(b,b) \\
 & \geq & f\circ g^*
\end{eqnarray*}  
as wanted.
\endofproof

In the next proposition, $\Matr(\Q)$ denotes the quantaloid of $\Q$-typed sets and matrices with elements in $\Q$ between them. We write a matrix typically as $M\:X\to Y$ (and its elements are $\Q$-arrows $M(y,x)\:tx\to ty$), and its composition with another matrix $N\:Y\to Z$ as $N\circ N\:X\to Z$ (with elements $(N\circ M)(z,x)=\bigvee_{y\in Y}N(z,y)\circ M(y,x)$). (Thus, $\Matr(\Q)$ is precisely the quantaloid of discrete $\Q$-enriched categories and distributors between them.)
\begin{proposition}\label{8.2}
A small involutive quantaloid $\Q$ is locally localic and modular if and only if $\Matr(\Q)$ is modular.
\end{proposition}
\proof 
First suppose that $\Q$ is locally localic and modular, and let $M\:X\to Y$, $N\:Y\to Z$ and $P\:X\to Z$ be arrows in $\Matr(\Q)$: we must prove that, for any $x\in X$ and $z\in Z$,
$$\Big(\bigvee_{y\in Y}N(z,y)\circ M(y,x)\Big)\wedge P(z,x)\leq\bigvee_{y'\in Y}\Big[N(z,y')\circ\Big(M(y',x)\wedge\bigvee_{z'\in Z}N\o(y',z')\circ P(z',x)\Big)\Big].$$
By distributivity of $\wedge$ over $\bigvee$ in $\Q$, and modularity of $\Q$, this can straightforwardly be verified:
\begin{eqnarray*}
LHS
 & = & \bigvee_{y\in Y}\Big(N(z,y)\circ M(y,x)\wedge P(z,x)\Big) \\
 & \leq & \bigvee_{y\in Y}\Big[N(z,y)\circ\Big(M(y,x)\wedge N(z,y)\o\circ P(z,x)\Big)\Big] \\
 & = & \bigvee_{y\in Y}\Big[N(z,y)\circ\Big(M(y,x)\wedge N\o(y,z)\circ P(z,x)\Big)\Big] \\
 & \leq & \bigvee_{y\in Y}\Big[N(z,y)\circ\Big(M(y,x)\wedge\Big(\bigvee_{z'\in Z}N\o(y,z')\circ P(z',x)\Big)\Big)\Big] \\
 & = & RHS.
\end{eqnarray*}

Secondly, suppose that $\Matr(\Q)$ is modular. Certainly $\Q$ is modular too, for it is a full subcategory. To see that $\Q$ is locally localic, let $f,(g_i)_{i\in I}\in\Q(X,Y)$; we need to show that $f\wedge(\bigvee_ig_i)\leq\bigvee_i(f\wedge g_i)$ (the reverse inequality is trivial). To see this, we consider the following sets and matrices: 
\begin{itemize}
\item $\{X\}$, the singleton whose single element is of type $X$,
\item $\{Y\}$, the singleton whose single element is of type $Y$,
\item $I$, the index-set for which we set the type of each $i\in I$ to $Y$,
\item $F\:\{X\}\to\{Y\}$, the matrix whose single entry is $F(Y,X)=f$,
\item $G\:\{X\}\to I$, the matrix whose $i$th entry is $G(i,X)=g_i$,
\item ${\bf 1}\:I\to \{Y\}$, the matrix whose $i$th entry is ${\bf 1}(Y,i)=1_Y$.
\end{itemize}
By these definitions, the morphism $f\wedge(\bigvee_i g_i)\:X\to Y$ in $\Q$ is the single element of the $\Q$-matrix $F\wedge{\bf 1}\circ G\:\{X\}\to\{Y\}$. By the hypothetical modularity of $\Matr(\Q)$, the latter is less than or equal to ${\bf 1}\circ({\bf 1}\o\circ F\wedge G)$, whose single element in turn is $\bigvee_i(f\wedge g_i)$. 
\endofproof

Finally, we state a proposition concerning modularity and tabularity, which is proved with calculations in the style of [Freyd and Scedrov, 1990, pages 223--224]; this will be useful in Example \ref{6.4}.
\begin{proposition}\label{3.02}
If an involutive quantaloid $\Q$ is modular and tabular then it is locally localic.
\end{proposition}
\proof
First remark that, for any morphism $q\:X\to Y$ in $\Q$, the modular law implies
$$q=q\wedge q\leq q(1_X\wedge q\o q)\leq qq\o q.$$
(The condition that $q\leq qq\o q$ for any morphism $q$ in $\Q$, is sufficient for many applications of modular quantaloids; such a quantaloid $\Q$ is sometimes said to be ``weakly Gelfand''.) In particular is any endomorphism $m\:X\to X$ such that $m\leq 1_X$, necessarily an idempotent: $mm\leq m$ holds in general, and $m\leq mm\o m\leq mm$ follows from the argument above. It is then straightforward that the sublattice $\down 1_X\subseteq\Q(X,X)$ of endomorphisms on $X$ below $1_X$, is a locale: because $mn=m\wedge n$ for any $m,n$ below $1_X$.

Now, for two objects $X,Y\in\Q$, taking advantage of Lemma \ref{3} we can choose $f\:U\to X$ and $g\:U\to Y$ such that $f\o f\wedge g\o g=1_U$, $fg\o=\top_{Y,X}$, $f\dashv f\o$ and $g\dashv g\o$. Thus there are adjoint order-preserving functions
$$\down 1_U\xymatrix@=12ex{\ar@/^2ex/[r]^{m\mapsto m}\ar@{}[r]|{\bot} & \ar@/^2ex/[l]^{1_U\wedge e\otspam e} }\Q(U,U)\xymatrix@=12ex{\ar@/^2ex/[r]^{e\mapsto feg\o}\ar@{}[r]|{\bot} & \ar@/^2ex/[l]^{f\o q g\otspam q}}\Q(Y,X)$$
so if we prove that the unit and counit inequalities of the composed adjunction 
$$\down 1_U\xymatrix@=12ex{\ar@/^2ex/[r]^{m\mapsto fmg\o}\ar@{}[r]|{\bot} & \ar@/^2ex/[l]^{1_U\wedge f\o qg\otspam q} }\Q(Y,X)$$
are in fact equalities, then $\Q(Y,X)$ is isomorphic (qua ordered set, thus also qua lattice) to the locale $\down 1_U$, and hence itself a locale. The inverse of the counit is easy to check: using modularity twice, we have
$$q=T_{Y,X}\wedge q=fg\o\wedge q\leq f(g\o\wedge f\o q)\leq f(1_U\wedge f\o qg)g\o.$$
For the inverse of the unit, first observe that
$$1_U\wedge f\o fmg\o g\leq 1_U\wedge f\o fmm\o g\o g=1_U\wedge (f\o fm)(g\o gm)\o$$
because $m\o\leq(1_U)\o=1_U$. But for any morphisms $a,b\in\Q(V,U)$ we can compute with the modular law that $1_U\wedge ab\o=1_U\wedge(1_U\wedge ab\o))\leq 1_U\wedge a(a\o\wedge b\o))\leq(a\wedge b)(a\wedge b)\o$. In our situation this implies that 
$$1_U\wedge (f\o fm)(g\o gm)\o\leq(f\o fm\wedge g\o gm)(f\o fm\wedge g\o gm)\o.$$ 
But because $(f\o fm\wedge g\o gm)\leq(f\o f\wedge g\o gmm\o)m\leq(f\o f\wedge g\o g)m=1_Um=m\leq 1_U$ (using the modular law, the fact that $mm\o\leq 1_U$, and the hypotheses on $f$ and $g$) we find 
$$(f\o fm\wedge g\o gm)(f\o fm\wedge g\o gm)\o\leq mm\o\leq m$$ 
as needed to conclude.
\endofproof

\section{Elementary characterisation of $\R(\C,J)$}\label{D}

If $\Q$ is small then we can regard $\Map(\Q)$ as a category and construct the quantaloid $\R(\Map(\Q))$ of {\em cribles of left adjoints} in $\Q$. To compare $\R(\Map(\Q))$ with the given $\Q$, there is always the normal colax functor $F\:\R(\Map(\Q))\to\Q$ defined to send a crible of left adjoints in $\Q$, say $R\:D\crible C$, to the $\Q$-morphism
$$F(R):=\bigvee\left\{fg^*\bigmid (f,g)\in R\right\}\in\Q(D,C).$$
For any objects $X$ and $Y$,
$$\R(\Map(\Q))(Y,X)\to\Q(Y,X)\:R\mapsto F(R)$$
preserves arbitrary suprema, hence admits a right adjoint qua order-preserving function: 
$$\Q(Y,X)\to\R(\Map(\Q))(Y,X)\:q\mapsto F^*(q).$$ 
Explicitly,
$$F^*(q)=\left\{(f,g)\in\R(\Map(\Q))(Y,X)\bigmid fg^*\leq q\right\}$$ 
and it follows easily that this defines a lax functor $F^*\:\Q\to\R(\Map(\Q))$. In the next two lemmas (the first of which is a mere triviality) we establish a link with the conditions in Definition \ref{2}.
\begin{lemma}\label{4.0} For a small quantaloid $\Q$, the following are equivalent:
\begin{enumerate}
\item $\Q$ is weakly tabular,
\item for all objects $X,Y$ in $\Q$, the adjunction
$$\R(\Map(\Q))(Y,X)\xymatrix@=7ex{\ar@/^2ex/[r]^{F}\ar@{}[r]|{\perp} & \ar@/^2ex/[l]^{F^*}}\Q(Y,X)$$
is split (ie.\ $F(F^*(q))=q$ for all $q\:Y\to X$ in $\Q$), 
\item $F\:\R(\Map(\Q))\to\Q$ is full.
\end{enumerate}
\end{lemma}
\begin{lemma}\label{4}
For a small weakly tabular and map-discrete quantaloid $\Q$, $F\:\R(\Map(\Q))\to\Q$ is a $\Sup$-functor and
$$j\:\R(\Map(\Q)\to\R(\Map(\Q))\:\Big(R\:Y\crible X\Big)\mapsto\Big(F^*(F(R))\:Y\crible X\Big)$$ 
is a nucleus such that $\Q\cong\R(\Map(\Q))_j$.
\end{lemma} 
\proof 
To prove that $F$ is functorial we must show that $F$ is lax on composites (for it is always normal colax): $F(R)\circ F(S)\leq F(R\circ S)$ for any $R\:Y\crible X$ and $S\:Z\crible Y$ in $\R(\Map(\Q))$. Equivalently: if $(f,g)\in R$ and $(m,n)\in S$ then $f g^* m n^*\leq F(R\circ S)$. Now, by weak tabularity of $\Q$ we know that $g^*m=\bigvee F^*(g^*m)$, so 
$$f g^* m n^*=\bigvee\left\{f a b^* n^*\bigmid (a,b)\in F^*(g^*m)\right\}.$$
But, for any $(a,b)\in F^*(g^*m)$,
$$ab^*\leq g^*m\ \Rightarrow \ ga\leq mb\ \Rightarrow\ ga=mb$$
since $\Q$ is map-discrete. From $(fa,ga)\in R$, $(mb,nb)\in S$ and $ga=mb$, it further follows that $(fa,nb)\in R\circ S$, whence $fab^*n^*\leq F(R\circ S)$. Thus we obtain $fg^*mn^*\leq F(R\circ S)$ as wanted. 

Secondly, $j=F^*\circ F$ is a nucleus on $\R(\Map(\Q))$ because it is a lax functor (it is the composite of two lax functors) and because locally, for any objects $X$ and $Y$,
$$\R(\Map(\Q))(Y,X)\to\R(\Map(\Q))(Y,X)\:R\mapsto j(R)$$ 
is a closure operator (it is the composite of the left and right adjoint in Lemma \ref{4.0}). By Lemma \ref{4.0} it is furthermore clear that the quotient quantaloid $\R_j(\Map(\Q))$ is isomorphic to $\Q$: the restriction of $F\:\R(\Map(\Q))\to\Q$ to the $j$-closed cribles is the identity on the objects, and fully faithful on the morphisms.
\endofproof
We can now prove our main result:
\begin{theorem}\label{6}
For a small quantaloid $\Q$, the following conditions are equivalent:
\begin{enumerate}
\item $\Q$ is locally localic, map-discrete, weakly tabular and weakly modular,
\item $F\:\R(\Map(\Q))\to\Q$ is a full and locally left exact $\Sup$-functor,
\item putting, for $X\in\Map(\Q)$,
$$J(X):=\Big\{S\mbox{ is a sieve on }X\bigmid 1_X=\bigvee_{s\in S}ss^*\Big\}$$
defines a Grothendieck topology $J$ on $\Map(\Q)$ for which $\Q\cong\R(\Map(\Q),J)$,
\item $\Q$ is a small quantaloid of closed cribles.
\end{enumerate}
In this case, $\Q$ carries an involution, sending $q\:Y\to X$ to
$$q\o:=\bigvee\left\{gf^*\bigmid (f,g)\:Y\span X\mbox{ is a span of left adjoints such that }fg^*\leq q\right\},$$
which makes $\Q$ a modular quantaloid.
\end{theorem} 
\proof
($1\Rightarrow 2$) Lemma \ref{4} provides everything except for the local left exactness of $F$. Thus it remains to prove that $F(R)\wedge F(S)=F(R\cap S)$ holds for morphisms $R,S\in\R(\Map(\Q))(Y,X)$. But $F(R)\wedge F(S)\geq F(R\cap S)$ is trivial (because $F$ is a $\Sup$-functor), and to check the other inequality it suffices -- because $\Q$ is locally localic -- to prove that $fg^*\wedge mn^*\leq F(R\cap S)$ for any $(f,g)\in R$ and $(m,n)\in S$. By weak modularity we have $f g^*\wedge m n^*\leq f(g^* n\wedge f^* m) n^*$, and by weak tabularity we further know that $g^* n\wedge f^* m=\bigvee F^*(g^* n\wedge f^* m)$, so it really suffices to prove that $fab^*n^*\leq F(R\cap S)$ for any $(a,b)\in F^*(g^* n\wedge f^* m)$. With an argument similar to that in the proof of Lemma \ref{4}, using in particular map-discretess, it is easily seen that $(a,b)\in F^*(g^* n\wedge f^* m)$ implies $f a=m b$ and $g a=n b$. Hence, from $(f a,g a)\in R$ and $(m b,n b)\in S$ it follows that $(f a,n b)\in R\cap S$, which implies $f a b^* n^*\leq F(R\cap S)$ as wanted.\\
($2\Rightarrow 3$) The nucleus $j\:\R(\Map(\Q))\to\R(\Map(\Q))$ of Lemma \ref{4} is locally left exact because it is the composite of locally left exact lax functors ($F^*$ even preserves all local infima). But $\Q\cong\R_j(\Map(\Q))$, as Lemma \ref{4} attests, and $\R_j(\Map(\Q))$ is necessarily the small quantaloid of closed cribles $\R(\Map(\Q),J)$ for the unique Grothendieck topology on $\Map(\Q)$ corresponding with the locally left exact nucleus $j$: thus
$$J(X)=\left\{S\mbox{ is a sieve on }X\bigmid\id_X\leq j(\{(s,s)\mid s\in S\})\right\}$$
which is precisely the same thing as in the statement of the theorem.\\
($3\Rightarrow 4$) Holds by Definition \ref{1}.\\
($4\Rightarrow 1$) Was explained in Example \ref{3.1}.\\
Finally, the involution $q\mapsto q\o$ results from the isomorphism in the third statement.
\endofproof

We can say a bit more about the topology constructed in the previous theorem.
\begin{proposition}\label{7.0}
If $\Q$ is a small quantaloid of closed cribles, then the Grothendieck topology $J$ on $\Map(\Q)$ as in Theorem \ref{6} is subcanonical (i.e.\ each representable is a sheaf).
\end{proposition}
\proof
Suppose that $S\in J(C)$ is a covering sieve, thus $\bigvee_{s\in S}ss^*=1_C$. With the usual abuse of notation we shall write $\sigma\:S\tto\Map(\Q)(-,C)$ for this sieve viewed as subfunctor of a representable functor, with $S(X)=\{s\in S\mid \dom(s)=X\}$ and $\sigma_X(s)=s$. For any other natural transformation into a representable, say $\tau\:S\tto\Map(\Q)(-,D)$, we must exhibit a unique morphism $f\:C\to D$ in $\Map(\Q)$ such that $\Map(\Q)(-,f)\circ\sigma=\tau$. The latter condition means precisely that $f\sigma_X(s)=\tau_X(s)$ for each $X\in\Map(\Q)$ and $s\in S(X)$. Keeping in mind that $\sigma_X(s)$ and $\tau_X(s)$ are left adjoints in $\Q$, it follows that
$$f=\bigvee_{X,s}\tau_X(s)\sigma_X(s)^*\mbox{ \ \ and  \ \ }f^*=\bigvee_{X,s}\sigma_X(s)\tau_X(s)^*$$
form the unique possible candidate for an adjunction $f\dashv f^*$ in $\Q$ satisfying the commutativity condition $f\sigma_X(s)=\tau_X(s)$. (In these suprema, $X$ ranges over all objets of $\Map(\Q)$ and $s$ ranges over all elements of $S(X)$. This notational convention reappears in the suprema below.) We complete the proof by checking that this $f$ does indeed meet these requirements:

First, the commutativity condition. In one direction we trivially have
$$\tau_X(s)\leq\tau_X(s)\sigma_X(s)^*\sigma_X(s)\leq\bigvee_{Y,t}\tau_Y(t)\sigma_Y(t)^*\sigma_X(s)=\Big(\bigvee_{Y,t}\tau_Y(t)\sigma_Y(t)^*\Big)\sigma_X(s)=f\sigma_X(s).$$
For the other direction, it suffices to show that $\tau_Y(t)\sigma_Y(t)^*\sigma_X(s)\leq\tau_X(s)$ for all $X,Y\in\Map(\Q)$ and $s\in S(X)$, $t\in S(Y)$, or equivalently $\sigma_Y(t)^*\sigma_X(s)\leq\tau_Y(t)^*\tau_X(s)$, or still equivalently, $t^*s\leq\tau_Y(t)^*\tau_X(s)$. We use the same trick as in the proofs of Lemma \ref{4} and Theorem \ref{6}: if $(a,b)\:X\span Y$ is a span of left adjoints in $\Q$ such that $ab^*\leq t^*s$, then $ta=sb$ holds by map-discreteness of $\Q$; and because $ta=sb$ is an element of the sieve $S$, we infer by naturality of $\tau$ that $\tau_Y(t)a=\tau_X(s)b$; this, in turn, is equivalent to $ab^*\leq\tau_Y(t)^*\tau_X(s)$. Because $\Q$ is weakly tabular, this suffices to prove that $t^*s\leq\tau_Y(t)^*\tau_X(s)$, as wanted.

Next, the unit of the adjunction. This is easy:
\begin{eqnarray*}
f^*f
 & = & \Big(\bigvee_{X,s}\sigma_X(s)\tau_X(s)^*\Big)\Big(\bigvee_{Y,t}\tau_Y(t)\sigma_Y(t)^*\Big) \\
 & \geq & \bigvee_{X,s}\sigma_X(s)\tau_X(s)^*\tau_X(s)\sigma_X(s)^* \\
 & \geq & \bigvee_{X,s}\sigma_X(s)\sigma_X(s)^* \\
 & = & \bigvee_sss^* \\
 & = & 1_C.
\end{eqnarray*}

Finally, the counit of the adjunction. We must show that $ff^*\leq 1_D$, that is, 
$$\Big(\bigvee_{X,s}\tau_X(s)\sigma_X(s)^*\Big)\Big(\bigvee_{Y,t}\sigma_Y(t)\tau_Y(t)^*\Big)\leq 1_D.$$ 
It suffices to show that $\tau_X(s)\sigma_X(s)^*\sigma_Y(t)\tau_Y(t)^*\leq 1_D$ for any $X,Y\in\Map(\Q)$ and $s\in S(X)$, $t\in S(Y)$, or equivalently $\sigma_X(s)^*\sigma_Y(t)\leq \tau_X(s)^*\tau_Y(t)$. But we have already shown this when checking the commutativity condition.
\endofproof
We can improve this result, provided that each {\em coreflexive arrow} $e\:C\to C$ (meaning that $e\leq 1_C$) in $\Q$ splits: there exist arrows $g\:E\to C$ and $f\:C\to E$ satisfying $fg=1_E$ and $gf=e$. In this case, it trivially follows that, necessarily, $g\dashv f$ and $e^2=e$. 
\begin{proposition}\label{7.01}
If $\Q$ is a small quantaloid in which coreflexive arrows split and $J'$ is a subcanonical Grothendieck topology on $\Map(\Q)$, then for each $C\in\Map(\Q)$ and each $S\in J'(C)$ we have $\bigvee_{s\in S}ss^*=1_C$. In particular, if $\Q$ is a small quantaloid of closed cribles in which coreflexive arrows split, then the Grothendieck topology $J$ on $\Map(\Q)$ as in Theorem \ref{6} is canonical.
\end{proposition}
\proof
For a covering sieve $S\in J'(C)$ (which, as before, we write as $\sigma\:S\tto\Map(\Q)(-,C)$ whenever this is useful), it is trivial that $e:=\bigvee_{s\in S}ss^*$ is a coreflexive arrow in $\Q$; thus, by assumption, it splits: there exists a left adjoint $g\:E\to C$ in $\Q$ such that $gg^*=e$ and $g^*g=1_E$. For each $s\in S(X)$ we have $es=s$, and it follows that
$$g^*ss^*g\leq g^*g=1_E\mbox{ \ \ and \ \ }s^*gg^*s=s^*es=s^*s\geq 1_X,$$
that is to say, $g^*s\dashv s^*g$. This provides for a natural transformation $\tau\:S\tto\Map(\Q)(-,E)$ with components $\tau_X(s):=g^*s$, and since it is easily checked that $gg^*s=es=s$ for each $s\in S$, we have $\Map(\Q)(-,g)\circ\tau=\sigma$. On the other hand, because $\sigma\:S\tto\Map(\Q)(-,C)$ is part of the subcanonical topology $J'$, the sheaf condition for the representable $\Map(\Q)(-,E)$ dictates the existence of a unique $f\:C\to E$ in $\Map(\Q)$ such that $\Map(\Q)(-,f)\circ\sigma=\tau$. The representable $\Map(\Q)(-,C)$ too is a sheaf for $J'$, so $\Map(\Q)(-,gf)\circ\sigma=\sigma$ implies that $gf=1_C$. Together with $g^*g=1_E$ this in turn means that $fg=g^*gfg=g^*g=1_E$, so that $f=g^{-1}=g^*$. We conclude that $1_C=gg^*=e$, as required.

The second part of the proposition is now an evident consequence of Proposition \ref{7.0}. 
\endofproof

To end this paper, we illustrate Theorem \ref{6} with some examples.
\begin{example}\label{6.4}
Suppose that $\Q$ is a small involutive quantaloid, with involution $q\mapsto q\obis$, and suppose that -- for this involution -- $\Q$ is modular and tabular. By Lemmas \ref{3.01}, \ref{3} and Proposition \ref{3.02} it follows that $\Q$ is locally localic, map-discrete, weakly modular and weakly tabular, i.e.\ a small quantaloid of closed cribles. But Theorem \ref{6} then implies that $\Q$ comes equipped with an involution $q\mapsto q\o$. It is easy to check that the two involutions on $\Q$ coincide: $q\o=q\obis$. 
\end{example}
\begin{example}\label{6.4.1}
The category $\Map(\Q)$ necessarily has binary products whenever $\Q$ is a modular and tabular small quantaloid. (Indeed, for objects $X,Y$ of $\Map(\Q)$, let $f\:P\to Y$ and $g\:P\to X$ form a tabulation in $\Q$ of the top element $\top_{Y,X}\in\Q(Y,X)$, then $(P,f,g)$ is the product of $X$ and $Y$ in $\Map(\Q)$: if $(Q,k,l)$ is another cone, then $f\o k\wedge g\o l$ is its unique factorisation through $(P,f,g)$.) However, this need not be the case if $\Q$ is a small quantaloid of closed cribles: consider for instance the category $\C$ whose only non-identity arrows are 
$$\xymatrix@C=4mm{
 & Z \\
X\ar[ru]^{f} & & Y\ar[lu]_{g}}$$
then, for the discrete topology $J$ on $\C$, the category $\Map(\R(\C,J))$ does not have the product $X\times Y$ (because there are no left adjoints in $\R(\C,J)=\R(\C)$ with codomains $X$ and $Y$). Hence $\R(\C,J)$ is a small quantaloid of closed cribles which is not tabular. 
\end{example}
This example shows that a small quantaloid of closed cribles is neither a {\em cartesian bicategory} [Carboni and Walters, 1987; Carboni {\it et al.}, 2008], nor a {\em tabular allegory} [Freyd and Scedrov, 1990], although it is a related concept. 
\begin{example}\label{6.2}
A locale $(L,\bigvee,\wedge,\top)$ (viewed as a monoid in $\Sup$, i.e.\ a one-object quantaloid) is not weakly tabular; it is thus not a quantale of closed cribles. However, its split-idempotent completion $L\si$ is a small quantaloid of closed cribles: all axioms are easy to verify. Furthermore, $\Map(L\si)\cong L$ (viewing the ordered set $L$ as category), and the Grothendieck topology constructed in Theorem \ref{6} is precisely the canonical topology associated to the locale $L$.
\end{example}
\begin{example}\label{6.3}
Let $\G$ be a small groupoid, and let $J$ be the smallest Grothendieck topology on $\G$: the small quantaloid of closed cribles $\R(\G,J)$ then equals the quantaloid of cribles $\R(\G)$. The latter in turn is isomorphic (as involutive quantaloid) to the free quantaloid $\Q(\G)$ on $\G$, equipped with its canonical involution $S\mapsto S\o:=\{s\inv\mid s\in S\}$. Indeed, any crible $R\:X\crible Y$ in $\G$ determines the subset $F(R):=\{h\inv g\mid (g,h)\in R\}$ of $\G(X,Y)$. Conversely, for any subset $S$ of $\G(X,Y)$ let $G(S)$ be the smallest crible containing the set of spans $\{(1_X,s)\mid s\in S\}$ in $\G$. Then $R\mapsto F(R)$ and $S\mapsto G(S)$ extend to functors $F\:\R(\G)\to\Q(\G)$ and $G\:\Q(\G)\to\R(\G)$ which are each other's inverse, and which preserve the respective involutions.
\end{example}
\begin{example}\label{6.1}
The quantale of extended positive real numbers $([0,\infty],\bigwedge,+,0)$ (viewed as a one-object quantaloid) is not weakly tabular; therefore it is not a quantale of closed cribles. As this quantale is equivalent to its split-idempotent completion $[0,\infty]\si$, the latter cannot be a small quantaloid of closed cribles either. 
\end{example}

\end{document}